\documentclass[11pt]{amsart}
\usepackage{amssymb,amscd}
\usepackage{epsf}

\newtheorem{thm}{Theorem}

\newtheorem{lemma}{Lemma}
\newtheorem{coro}{Corollary}

\begin{document}

\title{A knotted minimal tree}
\author{Krystyna Kuperberg}
\address{Auburn University, Auburn, AL 36830-5310, USA}
\email{kuperkm@math.auburn.edu}
\thanks{This research was supported in part by NSF grant \# DMS-9401408.}

\begin{abstract} There is a finite set of points on the boundary of the
three-dimensional unit ball whose minimal tree is knotted. This example answers
a problem posed by Michael Freedman. 
\end{abstract}

\maketitle

\section{Introduction}

\bigskip

In [1], ``Problems in Low-dimensional Topology'' by Rob Kirby,  one can find
the following:

\medskip 

\noindent{\bf Problem 5.17 (Freedman)} {\em  Given a finite set of points $X$
in $\partial B^3$, let $T$ be a tree in $B^3$ of minimal  length with
$\partial T=X$. Is $T$ unknotted, that is, is there a ${\rm PL}$  imbedded
2-ball in $B^3$ containing $T$?}

\medskip 

It is shown here that there is a finite set $X$ on the boundary unit 3-ball in ${\mathbb R}^3$ whose
minimal tree {\em is} knotted. The cardinality of $X$ is quite large, but the
construction essentially depends on seven elements only: 6 points and an arc on the
boundary of the ball.  An outline of the example, provided by the author of this paper, is contained
in [1] following the statement of Problem 5.17. 
 
\vspace{.2in}
\centerline{\epsfbox{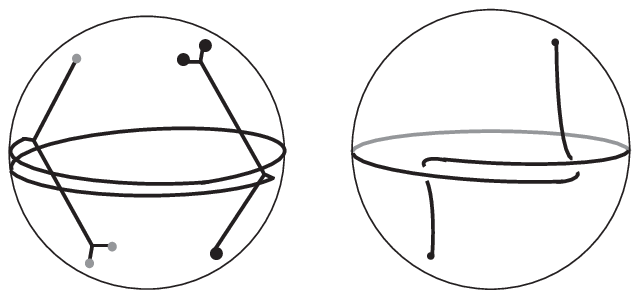}}
\vspace{.2in}

The arc lies close to the  equator circumventing it one and a half times.  In a
plane perpendicular to the equator, consider  a regular hexagon whose one pair of antipodal
vertices is very close to the endpoints of the arc. The set $X$ contains one  of the
other pairs of antipodal  vertices.  The remaining two vertices of the hexagon
are split into two  points each so that the points are closer to the equator and the
minimal tree connects these points to the endpoints of the
arc.  Finally, the arc is  replaced by a sequence of points of small mesh.

\section{Preliminaries and notation}

\bigskip 
Let $A$ be a locally connected compact set in ${\mathbb R}^3$ with finitely many
components. A {\em connecting graph\/} for $A$ is a pair $({\mathcal
E},{\mathcal V})$, where ${\mathcal E}$ is a finite collection of straight line
segments (edges) and ${\mathcal V}$ is a set of points (vertices) consisting of
endpoints of edges of ${\mathcal E}$, such that the set $A\cup
\bigcup_{E\in {\mathcal E}} E$ is connected. We informally say that the graph consists of
these edges. The {\em length\/} of a connecting graph is the sum of the lengths of its
edges. A  {\em minimal graph\/} for $A$ is a connecting graph for $A$ whose  length is a
minimum.  A minimal graph for
$A$ is denoted  by $G(A)$, possibly with a subscript if more than one of
such graphs is considered. The length of $G(A)$ is denoted by $|G(A)|$. No
additional notation is used for the union of the elements of $G(A)$ -- it is
also denoted by $G(A)$. For simplicity, assume that every point of $A\cap G(A)$
is a vertex of $G(A)$  but otherwise $G(A)$ has the minimum number of
vertices, i.e., if two edges meeting at a vertex are collinear, then the
vertex belongs to $A$. The {\em order\/} of a vertex is the number of edges
meeting at this vertex.  

If $A$ is finite, then a minimal graph $G(A)$ is a tree, i.e., it is connected
and acyclic. It is then called a {\em minimal tree\/} for $A$ and it is denoted 
by $T(A)$. Its length is denoted by $|T(A)|$. A {\em simple triod\/} is a tree
consisting of three edges meeting at a vertex. For a three point set
$A=\{a,b,c\}$, $T(A)$ is unique. If one of the angles of the triangle $\triangle
(a,b,c)$ is greater than or equal to
$\frac{2\pi}{3}$, then the minimal tree consists of two edges. If all angles of
$\triangle (a,b,c)$  are less than $\frac{2\pi}{3}$, then the minimal tree is  a
simple triod whose edges form
$\frac{2\pi}{3}$ angles. In general, a minimal tree is not  unique. For example,
the set of vertices of a square has two minimal trees.

If four or more half-lines in ${\mathbb R}^3$ have a common endpoint $p$,
then at least one of the angles between the half-lines is less than 
$\frac{2\pi}{3}$. If $q_1$ and $q_2$ are two distinct points equidistant to
$p$ that are on two half-lines meeting at $p$ at an angle less than 
$\frac{2\pi}{3}$, then the minimal tree $T(\{p,q_1,q_2\})$  is a triod.
Therefore, a vertex of a minimal graph $G(A)$ is either of order  3 or it
belongs to $A$.  The angles between the edges  meeting at a vertex not
in $A$ equal $\frac{2\pi}{3}$ and the edges are coplanar. 

The segment joining the points $p$ and $q$ is denoted by $[p,q]$.  The
Euclidean distance is denoted by ${\rm d}(p,q)$.  For distinct points $p$,
$q$ and $r$, denote by ${\rm L}(q,r)$ the line passing through $q$ and $r$,
and by ${\rm d}(p,{\rm L}(q,r))$ the  perpendicular distance between $p$ and
${\rm L}(q,r)$. The Hausdorff distance between the sets $A$ and $B$ is
denoted by ${\rm d_H}(A,B)$. We say that two sets $A$ and $B$ with the same
finite number of components, $A_1,\ldots ,A_k$ and $B_1,\ldots ,B_k$,
respectively,  are Hausdorff $\epsilon$-close, if there is a permutation
$\tau :\{1,\ldots ,k\}\to\{1,\ldots ,k\}$ such that for $i=1,\ldots k$,
${\rm d_H}(A_i,B_{\tau (i)})<\epsilon$. If $p$ and $q$ are non-antipodal
points on a circle or a sphere $C$, then ${\rm arc}_C(pq)$ denotes the shortest
arc in $C$ joining the points $p$ and $q$.

In here, the 3-ball $B^3$ is exactly the unit ball  in $ {\mathbb R}^3$. A PL imbedded
2-ball $D$ in $B^3$ is properly imbedded, i.e., $\partial B^3\cap D= \partial D$. A tree for
a finite set $A\subset\partial B^3$ is {\em unknotted\/} if there is a PL 2-ball  $D$
containing the tree, or equivalently, if there is an isotopy of $B^3$ onto itself such that
the image of the tree under the final stage of the isotopy is contained in the $xy$-plane.

Throughout the paper we use the following notation:
$$
\begin{array}{rcl}
S^2=\partial B^3&=&\{(x,y,z)\in {\mathbb R}^3\ |\ x^2+y^2+z^2=1\},\\
P&=&\{(x,y,z)\in S^2\ |\ y=0\},\\
Q&=&\{(x,y,z)\in S^2\ |\ z=0\},
\end{array}
$$
where $(x,y,z)$ denotes the Cartesian coordinates of a point in ${\mathbb R}^3$.

\section{Some special graphs}

Let $H\subset P$ be the regular hexagon with vertices:

$$
\begin{array}{lll}
a_1=(-\frac{1}{2},0 ,\frac{\sqrt{3}}{2}), 
&b_1=( -1,0,0),
&c_1=(-\frac{1}{2},0,-\frac{\sqrt{3}}{2}),\\
a_2=(\frac{1}{2},0,-\frac{\sqrt{3}}{2}),
&b_2=(1,0,0),
&c_2=(\frac{1}{2},0,\frac{\sqrt{3}}{2}). 
\end{array}
$$

An easy verification gives the following:

\begin{lemma}\label{four} Every  minimal tree $T(A)$  for the set $A$
consisting of four consecutive vertices of $H$ consists of three edges of $H$.
\end{lemma}

\begin{lemma}\label{five}    Every minimal tree $T(A)$ for the set $A$
consisting of five vertices of $H$ consists of four edges of $H$.
\end{lemma}

\begin{figure}[ht]
\centerline{\epsfbox{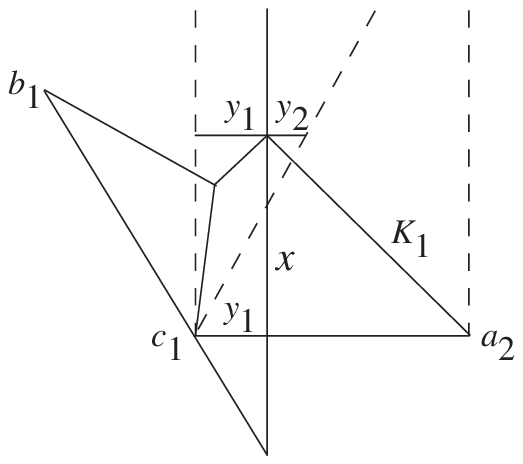}}
\caption{ $|T(A)|$ estimates.}
\label{ffivexy}
\end{figure}

\begin{proof} Let $A=\{a_1, b_1,c_1,a_2,c_2\}$. 

Suppose that $T(A)$ contains at least one edge of $H$. If $b_1$ is an order
1 vertex of $T(A)$ and it is one of the endpoints of the only edge of $H$
belonging to $T(A)$, then $|T(A)|=1+2\sqrt{3}$.  The cases when either $a_1$ or
$c_1$ is of order 1 can also be easily eliminated. If either $a_2$ or
$c_2$ is of order 1, then  Lemma ~\ref{four} can be used. 

\begin{figure}[ht]
\centerline{\epsfbox{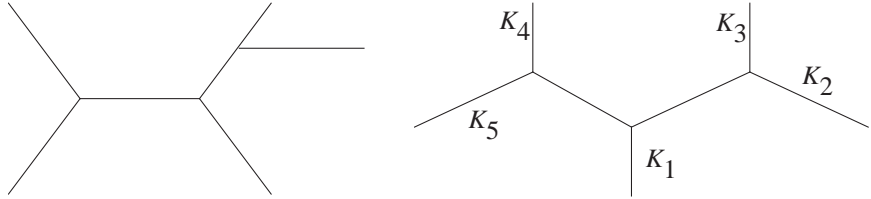}}
\caption{ The edges of $T(A)$ with three vertices of order 3.}
\label{finterior}
\end{figure}

Now suppose that no edge of $H$ belongs to $T(A)$, i.e., $T(A)$ has three additional
vertices. We may assume that every edge of $T(A)$ with an endpoint in
$H$ has length less than one. The combinatorial scheme for the edges of $T(A)$ is
as in Figure ~\ref{finterior}. The tree to the left shows that starting with any
edge joining two interior vertices and with the four adjacent edges, by attaching
the sixth edge  to one of the legs,  we obtain the tree to the right. The edge
$K_1$ is different from the other edges that have a vertex of $H$ as one of the endpoints:
$K_1$ does not have a common endpoint with any of the other edges $K_i$. Consider three
cases: $K_1=[p,b_1]$, $K_1=[p,c_1]$, and $K_1=[p,a_2]$ for some $p$ in the interior of the
hexagon, see Figures ~\ref{ffivebc} and ~\ref{ffivea}. Assume that all angles at interior
vertices equal $\frac{2\pi}{3}$.

Suppose that $K_1$ has $b_1$ as one of its endpoints. Then the pentagon whose
sides are the segment $[a_2,c_2]$ and four edges of $T(A)$, as shown in Figure
~\ref{ffivebc}, would have three $\frac{2\pi}{3}$ angles and the remaining two
angles each less than $\frac{\pi}{2}$. 

Suppose that $K_1$ has $c_1$ as one of its endpoints.  Then $T(A)$ contains two
paths, one from $a_1$ to $c_1$ and the other from $a_2$ to $c_2$  not overlapping
in a segment,  see Figure ~\ref{ffivebc}. Hence $|T(A)|> {\rm d}(a_1,a_2)+{\rm
d}(c_1,c_2)=4$.

\begin{figure}[ht]
\centerline{\epsfbox{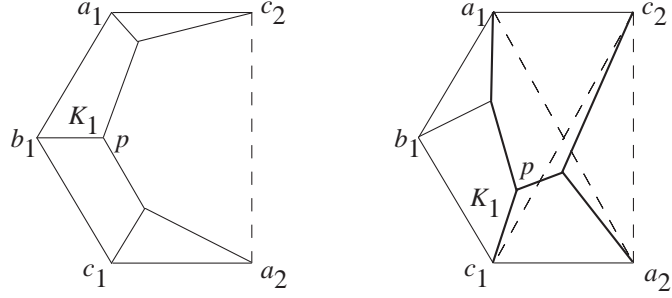}}
\caption{ Configuration $K_1=[p,b_1]$ and configuration $K_1=[p,c_1]$.}
\label{ffivebc}
\end{figure}

\begin{figure}[ht]
\centerline{\epsfbox{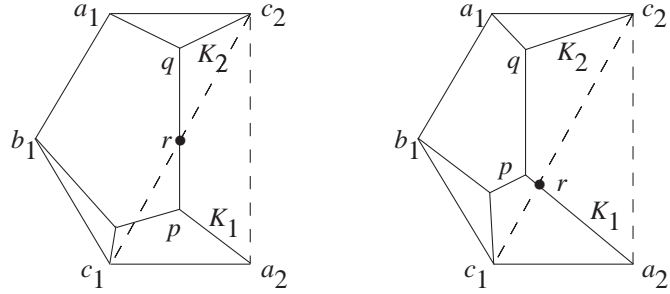}}
\caption{ Configuration  $K_1=[p,a_2]$.}
\label{ffivea}
\end{figure}

Suppose that $K_1$ has $a_2$ as one of its endpoints. Let $K_2=[q,c_2]$. The
angle $\angle(q,c_2,a_2)$ is greater than $\frac{\pi}{6}$, otherwise
$\angle(c_2,a_2,p)\geq\frac{\pi}{2}$.  The segment $[c_1,c_2]$ either
intersects the interior of $[p,q]$ or it intersects the edge $K_1$. In either
case denote the point of intersection by $r$.

If $r$ is in the interior of $[p,q]$, then by reflecting the part of $T(a)$ below
$r$ in the line passing through $c_1$ and $c_2$ we obtain a graph of the
same length as $T(a)$ and with an additional vertex at $r$, which is not
possible.

Suppose now that $[c_1,c_2]$ intersects $K_1$. Let $x={\rm d}(p,{\rm L}(c_1,a_2))$ and
$y={\rm d}(p_1,p_2)$, where the segment $[p_1,p_2]$ is parallel to ${\rm L}(c_1,a_2)$,
$p\in [p_1,p_2]$, $p_1\in {\rm L}(c_1,a_1)$, and $p_2\in {\rm L}(c_1,c_2)$. Let
$y_1={\rm d}(p_1,p)$ and $y_2={\rm d}(p_2,p)$, see Figure ~\ref{ffivexy}.

Note that 

\begin{enumerate}
\item since $|K_1|<1$, $p$ is on the same side of ${\rm L}(c_1,a_1)$ as $a_2$,
\item $|K_1|= {\rm d}(a_2,r)+{\rm d}(r,p)\geq \frac{\sqrt 3}{2}+ {\rm
d}(p,{\rm L}(c_1,c_2))\geq \frac{\sqrt 3}{2}+\frac{\sqrt 3}{2}y_2 $,
\item ${\rm d}(p,{\rm L}(b_1,c_1))=\frac{x+\sqrt{3}y_1}{2}$,
\item $y=\frac{x}{\sqrt 3}$,
\item $|T(\{p,a_1,c_2\})|\geq |G(\{a_1,c_2\}\cup{\rm L}(p_1,p_2))|=
\frac{3}{2}\sqrt{3}-x$,
\item $|T(\{p,b_1,c_1\})|\geq {\rm d}(p,{\rm L}(b_1,c_1))+\frac{\sqrt 3}{2} =
\frac{x}{2}+\frac{\sqrt 3}{2}y_1+\frac{\sqrt 3}{2}$\
\ (this holds true even if ${\rm d}(p,{\rm L}(b_1,c_1))<\frac{1}{2\sqrt 3}$ and
$T(\{p,b_1,c_1\}$ is not a triod).
\end{enumerate}
We have $ |T(A)|=|K_1|+ |T(\{p,a_1,c_2\})|+ |T(\{p,b_1,c_1\})| 
\geq \frac{\sqrt 3}{2} +\frac{\sqrt 3}{2}y_2 +
\frac{3}{2}\sqrt 3-x+ \frac{x}{2}+\frac{\sqrt 3}{2}y_1+ \frac{\sqrt 3}{2} =
\frac{5\sqrt 3}{2}- \frac{x}{2}+ \frac{\sqrt 3}{2}y= \frac{5\sqrt
3}{2}>4$.
\end{proof}

\begin{lemma}\label{hexa}   Every  minimal graph $G(A)$ for the set 
$A=Q\cup \{a_1,c_1,a_2,c_2\}$ consists of four edges of $H$.
\end{lemma}

\begin{proof} Let $b_4=(0,1,0)$ and  $b_5=(0,-1,0)$.

If $G(A)\cap Q\subset \{b_1, b_2\}$, then $G(A)$ is
planar; otherwise the projection of $G(A)$ onto the $xz$-plane would be of
shorter length. Then the vertices of a component of $G(A)$ that are in $H$ form a
sequence of consecutive vertices of $H$.   By the previous lemmas,
$G(A)$ consists of 4 edges of $H$.

Suppose that one of the edges $[a_1,c_2]$ or $[c_1,a_2]$  belongs to $G(A)$. The
only segments  in the convex hull of $A$ forming with these edges an angle
$\geq \frac{2\pi}{3}$ at one of the endpoints $a_1,a_2,c_1$, or $c_2$ are contained in
another edge of $H$. So if $G(A)$ contains either of these two edges, then it
contains an edge adjacent  to $[a_1,c_2]$ or $[c_1,a_2]$. If $G(A)$ does not
contain $[a_1,c_2]$ (resp. $[c_1,a_2]$) but it contains an edge  of $H$  adjacent
to it, then this edge can be replaced by  $[a_1,c_2]$ (resp. $[c_1,a_2]$) to get a
connecting graph of the same length.  Hence we may assume that if $G(A)$ contains an
edge of $H$, then it also contains another edge on the same side of the $xy$-plane.
Our consideration may be reduced to graphs $G(A)$ whose components contain either two
or four of the vertices $a_1,a_2,c_1,c_2$. The case when one component contains
$a_1$ and $a_2$, and the other $c_1$ and $c_2$ can be easily eliminated. If
one of  the components  contains $a_1$ and $c_1$, and  the other
contains  $a_2$ and $c_2$, then $G(A)$ is planar. Suppose that a non-planar
component of $G(A)$ contains exactly two vertices of $H$, both on the same side of
the $xy$-plane, say $a_1$ and $c_2$. Then $G(A)$ has a component that is a simple
triod with one additional endpoint $a_3\in Q$ different from $b_1$ and $b_2$, and a
vertex $s$ of order 3. Since  $a_3$ is the closest point to $s$ on $Q$, the line
${\rm L}(s,a_3)$ intersects the $z$-axis. The only possible choices for $a_3$ so
that $\angle (a_1,s,a_3)=\angle (c_2,s,a_3)$ are $b_4$ and $b_5$. But for $i=4,5$,
$|T(\{a_1,c_2,b_i\})|=\frac{\sqrt 7+\sqrt 3}{2}>2$.

Hence if $G(A)$ is non-planar, then $G(A)$ is connected and contains a point
$b_3\in  Q$ different from $b_1$ and $ b_2$. We may assume that $b_3\in  {\rm arc}_Q
(b_2,b_4)$.  $G(A)$ is combinatorially equivalent to the graph pictured in Figure
~\ref{finterior}. Let $p$ be the endpoint of $K_1$ that does not belong to $H$.
The possible types of configurations are:

\begin{enumerate}
\item $K_1=[b_3,p]$ and the endpoints of $K_2$ and $K_3$ that belong to $H$ are on
the same side of the $xy$-plane (see Figure ~\ref{fiveb}),
\item $K_1=[b_3,p]$ and the endpoints of $K_2$ and $K_3$ that belong to $H$ are on
the opposite sides of the $xy$-plane,
\item $K_1=[c_2,p]$, $K_2$  connects to $a_1$ and $K_3$ connects to $b_3$ (see
Figure ~\ref{fiveac}),
\item $K_1=[a_1,p]$, $K_2$  connects to $c_2$ and $K_3$ connects to $b_3$ (see
Figure ~\ref{fiveac}),
\item $K_1=[c_2,p]$, $K_2$  connects to $b_3$ and $K_3$ connects to either $c_1$ or
$a_2$,
\item $K_1=[a_1,p]$, $K_2$  connects to $b_3$ and $K_3$ connects to either $c_1$ or
$a_2$.
\end{enumerate}

\begin{figure}[ht]
\centerline{\epsfbox{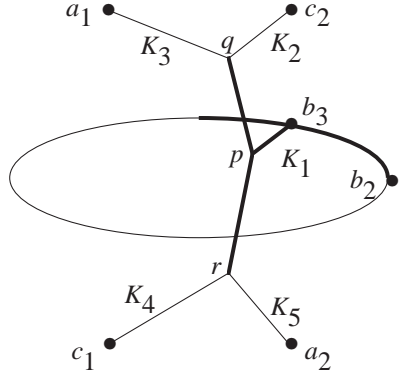}}
\caption{  $K_1$ connects to $b_3$.}
\label{fiveb}
\end{figure}

The most interesting is Configuration 1.  Let $q$ and $r$ be the remaining
two interior vertices different from $p$, with $q$ above the $xy$-plane and
$r$ below. Let $p',q',r',b'_3$ be the points obtained from $p,q,r,b_3$,
respectively, by a rotation of the tree $T(\{p,r,b_3\})$ in the $z$-axis so
that ${\rm d}(b'_3,b_2)<{\rm d} (b_3,b_2)$. Let $E$ be the ellipsoid given by
the equation ${\rm d}(x,a_1)+{\rm d}(x,c_2)=|K_2|+|K_3|$. Since $|K_2|={\rm
d}(c_2,q)\leq |K_3|={\rm d}(a_1,q)$, the point $q'$ is inside the ellipsoid
${E}$. Therefore ${\rm d}(a_1,q')+{\rm d}(c_2,q')<|K_2|+ |K_3|$. Similarly
${\rm d}(c_1,r')+{\rm d}(a_2,r')<|K_4|+ |K_5|$ and we obtain a connecting
graph for the set $A$ of shorter length. Hence $b_3=b_2$.

It is easy to eliminate Configurations 2, 5, and 6. In each of these cases 
$G(A)$ contains two non-overlapping paths from the set $\{a_1,c_2\}$ to the
set $\{a_2,c_1\}$, and since $G(A)$ is connected, $|G(A)|>4$.

\begin{figure}[ht]
\centerline{\epsfbox{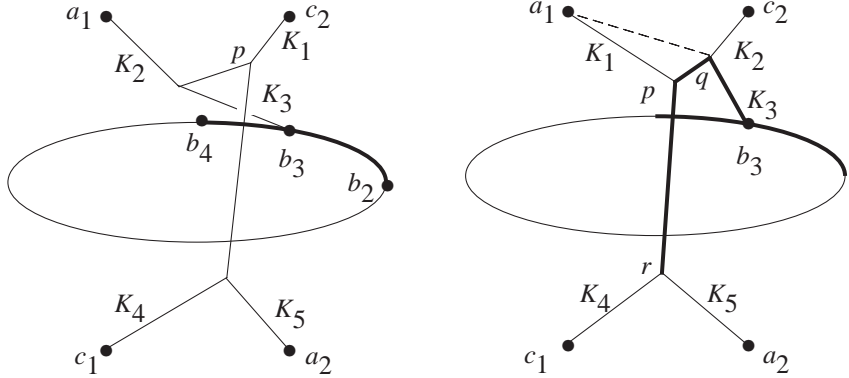}}
\caption{ $K_1$ connects to $c_2$; $K_1$ connects to $a_1$.}
\label{fiveac}
\end{figure}

Now consider Configuration 3.  We have,  
$$|G(A)|>|T(\{a_1,b_3,p\})| +|T(\{c_1,a_2,c_2\})|>$$
$${\rm d}(a_1,b_4)+|T(\{c_1,a_2,(0,0,\frac{\sqrt 3}{2})\})={\sqrt
2}+3\frac{\sqrt 3}{2}>4.$$

Finally consider  Configuration 4. Let $q$ and $r$ be the two additional interior
vertices with $K_2=[q,c_2]$ and $r$ being the common vertex of $K_4$ and $K_5$ as in
the right-hand picture in Figure ~\ref{fiveac}. Note that $b_3$ is the point
in $Q$ closest to $q$, hence the line $L(b_3,q)$ intersects the $z$-axis. If the
segment $[p,r]$ intersects the $yz$-plane at $\tilde r$ different from $r$, then by
reflecting the tree $T(\{ {\tilde r},a_2,c_1\})$ in the $yz$-plane, we obtain a
connecting graph for the set $A$ of the same length as $G(A)$ with an additional vertex
$\tilde r$. Therefore, we may assume that the points $p$ and $r$ are not on the opposite
sides of the $yz$-plane.

Suppose that $r$ is either on the $yz$-plane or  on the same side of the
$yz$-plane as $c_2$. Let $b'_3,p',q',r'$ be the new vertices corresponding to
$b_3,p,q,r$ obtained by  rotating the path $[b_3,q]\cup
[q,p]\cup [p,r]$ about the $z$-axis. The points $\overline a_1$,
$\overline c_2$, $\overline p$, $\overline p'$, $\overline q$, $\overline q'$ and
$\overline r$ in Figure ~\ref{frotate} are the projections of $a_1$, $c_2$, $p$,
$p'$, $q$, $q'$ and $r$ onto the $xy$-plane;  $o$ is the origin. Since $p$ is in the
convex hull of $\{a_1,c_1,a_2,c_2,b_3\}$ and $q$ is in the triangle $\triangle
(p,b_3,c_2)$,  then $\overline p\in \triangle (\overline a_1,o,\overline q)$. 

\begin{figure}[ht]
\centerline{\epsfbox{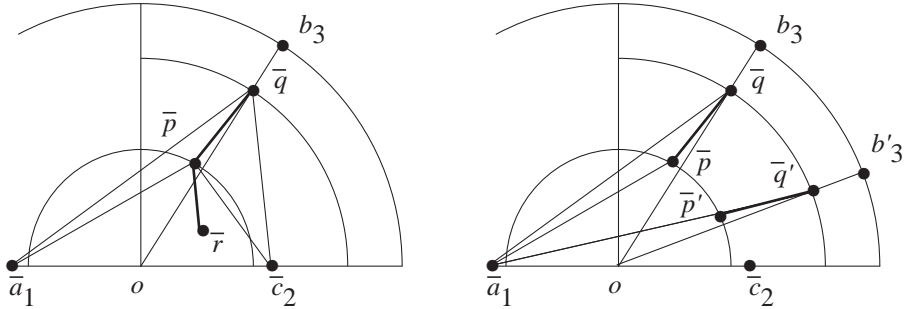}}
\caption{  Rotating part of $G(A)$.}
\label{frotate}
\end{figure}

If $b'_3\in {\rm arc}_Q(b_2,b_3)$ and $b'_3\not= b_3$, then,  arguing in a  similar
fashion as for Configuration 1, we get ${\rm d}(c_1,r')+{\rm d}(a_2,r')<
|K_4|+|K_5|$ and ${\rm d}(a_1,q')+{\rm d}(c_2,q') <{\rm d}(a_1,q)+{\rm d}(c_2,q)$.
Also note that ${\rm d}(a_1,q')>{\rm d}(a_1,q)$ and ${\rm d}(c_2,q')<{\rm
d}(c_2,q)$. Suppose that $\overline q$, $\overline p$ and $o$ are not collinear, and
$\overline q'$, $\overline p'$ and $\overline a_1$ are collinear as in the  right-hand
picture in Figure ~\ref{frotate}. Then $\angle (q',a_1,p')<\angle (q,a_1,p)$ and $\angle
(a_1,q',p')<\angle (a_1,q,p)$, and since ${\rm d}(p',q')={\rm d}(p,q)$, then ${\rm
d}(a_1,q')-{\rm d}(a_1,p')> {\rm d}(a_1,q)-{\rm d}(a_1,p)$, see Figure ~\ref{fangles}. 
Hence ${\rm d}(a_1,q')-{\rm d}(a_1,q)> {\rm d}(a_1,p')-{\rm d}(a_1,p)$, which combined with 
${\rm d}(a_1,q')+{\rm d}(c_2,q') <{\rm d}(a_1,q)+{\rm d}(c_2,q)$ gives ${\rm d}(a_1,p')+{\rm
d}(c_2,q') <{\rm d}(a_1,p)+{\rm d}(c_2,q) = |K_1|+|K_2|$.  We obtain a new connecting graph
for the set $A$ with interior vertices $p',q',r'$ whose length is less than $|G(A)|$.

\begin{figure}[ht]
\centerline{\epsfbox{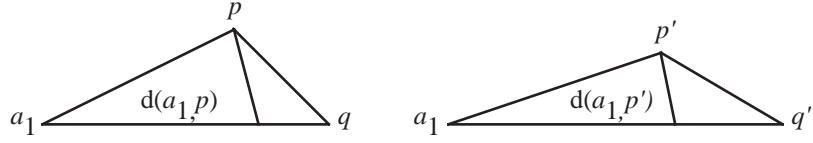}}
\caption{ ${\rm d}(a_1,q')-{\rm d}(a_1,p')> {\rm d}(a_1,q)-{\rm d}(a_1,p)$.}
\label{fangles}
\end{figure}

If the points $\overline a_1$, $\overline p$ and $\overline q$ are collinear,
then $a_1$, $p$ and $q$ are in a plane perpendicular to the $xy$-plane and so is the
tree $T(\{a_1,q,r\})$, in particular, so is the segment $[p,r]$. Then $G(A)$ is
planar.

If the points $\overline o$, $\overline p$ and $\overline q$ are collinear, then 
$b_3$, $q$, $p$ and $a_1$ are in a plane that is perpendicular to the $xy$-plane and
passes through the $z$-axis. This plane also contains $c_2$ and 
$G(A)$ is planar.

\begin{figure}[ht]
\centerline{\epsfbox{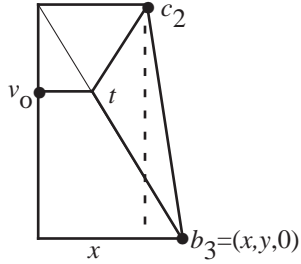}}
\caption{ $|G(A)|$ estimates.}
\label{ftrapezoid}
\end{figure}

If  $p$ and $r$ are on the same side or on the $yz$-plane as $a_1$ (including the
case when one or both of these points are on the $yz$-plane), then
$$|G(A)|>|T(\{a_1,c_1,a_2\})|+ {\rm min}_{v\in yz{\rm -plane}}|T(\{v,c_2,b_3\})|.$$
If ${\rm min}_{v\in yz{\rm -plane}}|T(\{v,c_2,b_3\})|$ is attained at $v_0\in
yz$-plane, then either
\begin{enumerate}
\item $v_0=b_3=b_4$, $|T(\{v_0,c_2,b_3\})|=\sqrt 2$ and $|G(A)|= \sqrt
2+\frac{3\sqrt 3}{2}>4$, or
\item  $v_0=(0,0,\frac{\sqrt 3}{2})$, $|T(\{v_0,c_2,b_3\})|=\frac{3}{2}$ and 
$|G(A)|= \frac{3}{2}+\frac{3\sqrt 3}{2}>4$ , or
\item   $|T(\{v_0,c_2,b_3\})|$ is a simple triod.
\end{enumerate}
In the last case, let $t$ be the order 3 vertex of $T(\{v_0,c_2,b_3\})$. Note that
the line ${\rm L}(v_0,t)$ is perpendicular to the $yz$-plane and the line ${\rm
L}(b_3,t)$ intersects the $z$-axis. The plane $K$ containing  ${\rm L}(v_0,t)$ and
${\rm L}(b_3,t)$ also contains the tree $T(\{v_0,c_2,b_3\})$, see Figure
~\ref{ftrapezoid}. If $b_3=(x,y,0)$, then
$\sqrt 3x=\sqrt{(\frac{\sqrt3}{2})^2+y^2}=\sqrt{\frac{7}{4}-x^2}$. Then
$x=\frac{\sqrt 7}{4}$,  $|T(\{v_0,c_2,b_3\})|=\frac{1}{4} +\frac{\sqrt 7}{2}$ and
$|G(A)|> 2(x- \frac{1}{2})+ \frac{1}{2}+ \frac{\sqrt
3}{2}(\sqrt{\frac{7}{4} -x^2}-\sqrt 3(x-\frac{1}{2})) =\frac{3\sqrt 3}{2}+
\frac{1}{4}+\frac{\sqrt 7}{2}>4$. 
\end{proof}

\begin{coro}\label{Q+graph} For every $\epsilon >0$ there is a $\delta >0$ such
that if  a set $A$ is Hausdorff $\delta$-close to $Q\cup \{a_1,c_1,a_2,c_2\}$, then
every minimal graph $G(A)$ is in an $\epsilon$-neighborhood of four edges of $H$. 
\end{coro}

In the following lemmas, the notation $Q_\delta$ is used for a circle in $B^3$
that is in a plane parallel to the $xy$-plane, with center on the $z$-axis, and
such that ${\rm d_H}(Q,Q_\delta )<\delta $.

\begin{lemma}\label{x} There are an $\epsilon >0$, a $\delta >0$, and 
an $\eta >0$  such that if  
\begin{enumerate}
\item $a'_1$ and $c'_2$ are points in the $xz$-plane such that ${\rm
d}(a_1,a'_1)<\eta$, 
${\rm d}(c_2,c'_2)<\eta$, and 
\item  $G(A)$ is a minimal graph for  $A=Q_\delta \cup 
\{a'_1,c'_2\}$ contained in an $\epsilon$-neighborhood
of the set $[a_1, b_1]\cup  [a_1,c_2]$,
\end{enumerate}
then $G(A)$ is contained in the  $xz$-plane.
\end{lemma}

\begin{proof} Let $q$ be the point of $G(A)$ that belongs to
$Q_\delta$. Thus $G(A)$ is the tree $T(\{q,a'_1,c'_2\})$. For  small
$\epsilon$, $\delta $, and $\eta$, the tree $T(\{q,a'_1,c'_2)$ consists of 2
segments in the $xz$-plane or it is a simple triod with an additional vertex $p$ 
close  to $a'_1$. Then the point of intersection $r$ of the line ${\rm L}(p,q)$ and the
edge $[a'_1,c'_1]$ is also close to $a'_1$. Since $q$ is the point on $Q_\delta$
that is closest to $p$,  the line ${\rm L}( p,q)$ intersects the $z$-axis at some
point $s$. Hence ${\rm L}(p,q)$ has two distinct points $r$ and $s$ in the
$xz$-plane and $T(\{q,a'_1,c'_2\})$ is in the $xz$-plane.
\end{proof}

For $0<\gamma <1$, let 
$c_1(\gamma)=(-\frac{1}{2},0,-\frac{\sqrt{3}}{2}(1-\gamma ))$ and 
$c_2(\gamma)=(\frac{1}{2},0,\frac{\sqrt{3}}{2}(1-\gamma ))$. 

\begin{lemma}\label{1-x} There is a
$\gamma_0 >0$,  such that the set $A=Q\cup  \{a_1,c_1(\gamma),a_2,c_2(\gamma)\}$,
${0<\gamma<\gamma_0}$, has a unique minimal graph $G(A)$ consisting of two
simple triods $$T(\{a_1,b_1,c_1(\gamma )\})\  and\ T(\{a_2,b_2,c_2(\gamma )\}).$$
\end{lemma}

\begin{proof} By Corollary ~\ref{Q+graph}, there are six cases of a minimal
graph $G(A)$ to consider.  $G(A)$ may be close to one of the sets consisting of
the following fours edges of $H$:

$$
\begin{array}{lll}
\vspace{.1in}
1.\  [a_1,b_1]\cup  [b_1,c_1]\cup  [a_2,b_2]\cup  [b_2,c_2],
&2.\  [a_1,b_1]\cup  [a_1,c_2]\cup  [c_1,a_2]\cup  [a_2,b_2],\\
\vspace{.1in}
3.\  [a_1,c_2]\cup  [c_2,b_2]\cup  [a_2,c_1]\cup  [c_1,b_1],
&4.\  [a_1,b_1]\cup  [b_1,c_1]\cup [a_1,c_2]\cup  [a_2,b_2],\\
\vspace{.1in}
5.\  [a_1,b_1]\cup  [b_1,c_1]\cup  [c_1,a_2]\cup  [b_2,c_2],
&6.\  [a_1,b_1]\cup  [b_1,c_1]\cup  [a_1,c_2]\cup  [c_1,a_2]. 
\end{array}
$$

The remaining cases of subsets of $H$ consisting of four edges are three cases
symmetric with respect to the origin to Cases 4, 5, or 6, and six possibilities of
graphs that do not connect to at least one of the vertices of $H$.

\begin{figure}[ht]
\centerline{\epsfbox{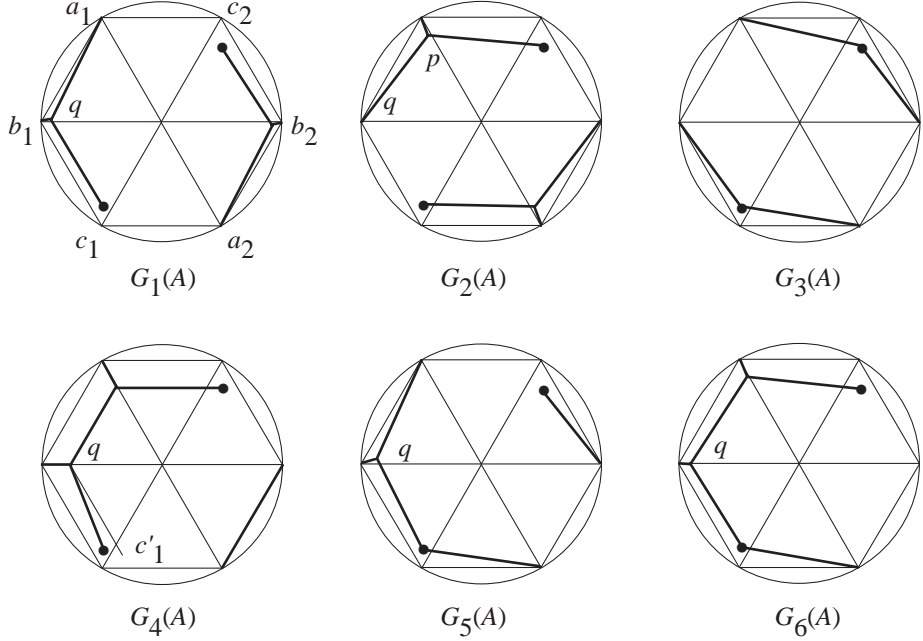}}
\caption{ The graphs $G_i(A)$.}
\label{QH}
\end{figure}

For $i=1,\ldots ,6$, let $G_i(A)$ be a minimal graph for $A$ corresponding to Case
$i$, see Figure  ~\ref{QH}.  Note that $G_1(A)$ and $G_3(A)$ are clearly contained
in the $xz$-plane. The see that the remaining graphs are also 
in the $xz$-plane, denote by  $q$ the additional vertex of the graph
$G_i(A)$ so that there is an edge  $[q,b_1]$ for $i=4,5,6$,   and
$q=b_1$ for $i=2$.  The vertex $q$ separates $G_i(A)$ into subgraphs. The two subgraphs 
different from the edge $[q,b_1]$ connect one or two of the points $a_1,c_1(\gamma),
a_2,c_2(\gamma)$ to a circle $Q_\delta$. By  Lemma ~\ref{x}, we may
assume that the subgraphs are subsets of the $xz$-plane. Since the distance between a
point in $Q_\delta$ and the circle $Q$ is constant, it easily follows that each $G_i(A)$ is
in the $xz$-plane. We analyze the six cases as follows:

\begin{enumerate}
\item $|G_1(A)|<2({\rm d}(a_1,b_1) + {\rm d}(b_1,c_1(\gamma)))< 2(1+1-\gamma +\frac{\gamma
}{2})=4-\gamma$.

\item $|G_2(A)| = 2({\rm d}(a_1,p)+{\rm d}(p,b_1)+{\rm d}(c_2(\gamma),p))>2(1+{\rm
d}(c_2(\gamma),p))$. Since  $p$ is on the same side of the
line ${\rm L}(a_1,a_2)$ as $q$, and ${\rm d}(p,{\rm L}(a_1,c_2))< {\rm d}(c_2(\gamma ),{\rm
L}(a_1,c_2))$, then 
$|G_2(A)|>2(1+1-\frac{\gamma }{2})>|G_1(A)|$.

\item $|G_3(A)|> 2({\rm d}(a_1,c_2(\gamma))+{\rm d}(c_2(\gamma),b_2))>2(1+{\rm
d}(c_2(\gamma),b_2))>|G_1(A)|$. 

\item  Let  $c'_1=(-\frac{1}{2}+\frac{3\gamma}{4},0,- \frac{\sqrt 3}{2}+\frac{\sqrt
3\gamma}{4})$, i.e.,
$c'_1$ is the point symmetric to
$c_1(\gamma)$ with respect to the line ${\rm L}(c_1,c_2)$. Let $C$ be the circle of radius
$\frac{\sqrt 3}{2}\gamma $ centered at $c'_1$. 
Since $c_1(\gamma)\in C$, $|G_4(A)|> |G(\{c_2(\gamma),a_1,b_1\}\cup  C)|-{\rm d}(c_1,c'_1)
+1 = 1-\frac{\gamma }{2}+1-\gamma +1-\frac{\gamma }{2}+ 2\gamma -\frac{\sqrt 3}{2}\gamma
+1=4-\frac{\sqrt 3}{2}\gamma >|G_1(A)|$.

\item   Suppose that $G_5(A)$ has two vertices $q$ and $u$ of order 3 and an edge
$[c_2(\gamma),u]$. Then the sum of the angles of the trapezoid with vertices
$a_1,q,u,a_2$ is less than $2\pi$. Hence $c_1(\gamma)$ is of order 2. We have
$|G_5(A)|=\frac{1}{2}|G_1(A)| +{\rm d}(c_1(\gamma),a_2) +{\rm
d}(c_2(\gamma),b_2)=\frac{1}{2}|G_1(A)|+\frac{1}{2}|G_3(A)|>|G_1(A)|$. 

\item   Suppose that $G_6(A)$ has three vertices $p$, $q$ and $u$ of order 3 and edges
$[a_1,p]$, $[b_1,q]$ and $[c_1(\gamma),u]$. The pentagon $U$ with vertices
$a_2,u,q,p,c_2(\gamma)$ is inside the  pentagon with vertices
$a_2,c_1(\gamma),b_1,a_1,c_2(\gamma)$. Since the edge $[a_2, c_1(\gamma)]$ is parallel
to the edge $[a_1, c_2(\gamma)]$, we have  $\angle (a_1, c_2(\gamma),a_2)+ \angle
(c_2(\gamma),a_2, c_1(\gamma))=\pi$ and  $\angle (p, c_2(\gamma),a_2)+ \angle
(c_2(\gamma),a_2, u)<\pi$. This is a contradiction since each of the remaining angles of
$U$ equals $\frac{2\pi}{3}$. Hence $c_1(\gamma)$ is of order 2. We have
$|G_6(A)|=|G_4(A)|-1 + {\rm d}(c_1(\gamma), a_2)>|G_4(A)|$. 

\end{enumerate}
\noindent Hence $G_1(A)$ is the unique minimal graph.

\end{proof}

For a given $\gamma$ and $i=1,2$,  let $L_i(\gamma )$ be the line passing through
$c_i(\gamma )$ and perpendicular to the
$xz$-plane. Denote by $e_i(\gamma )$ and $f_i(\gamma )$ the two
points  in which $L_i(\gamma )$ intersects the sphere $S^2$.

\begin{figure}[ht]
\centerline{\epsfbox{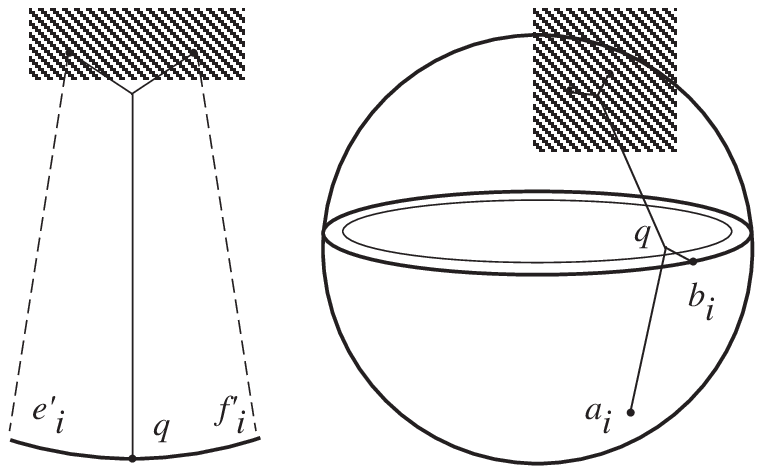}}
\caption{ Splitting the $c_i's$.}
\label{fsplit}
\end{figure} 

\begin{lemma}\label{split}  There are a $\gamma_0 >0$ and a $\delta >0$,  such that
for
$0<\gamma <\gamma_0$ and $i=1,2$, 
\begin {enumerate}
\item the minimal graph for the set $Q_\delta\cup \{e_i(\gamma ),f_i(\gamma )\}$ is the triod
$$T(\{v_i,e_i(\gamma),f_i(\gamma ) \}),$$
\item the minimal graph for the set $Q_\delta\cup \{a_i, e_i(\gamma ),f_i(\gamma )\}$ is the
minimal tree $$T(\{a_i, v_i, e_i(\gamma),f_i(\gamma )\}),$$
\end{enumerate}
where $v_i\in Q_{\delta}$ is  a point in the $xz$-plane close to $b_i$.
\end{lemma}
\begin{proof} For a given $\delta$, let $e'_i$ and $f'_i$ be points on a $Q_\delta$
closest to $e_i(\gamma )$ and $f_i(\gamma )$ respectively. Let $p\in {\rm arc}_{Q_\delta}
(e'_i,f'_i)$. The minimal tree $T(\{q,e_i(\gamma),f_i(\gamma ) \})$ has an edge containing $q$,
whose extension intersects the $z$-axis, see Figure ~\ref{fsplit}. For small
$\gamma$ and $\delta$ this is possible only for $q\in xz$-plane. From this, both
Conclusions 1 and 2 of Lemma ~\ref{split} follow.
\end{proof}

Similarly, one can prove the following:

\begin{coro}\label{splitmore}  There are a $\gamma_0 >0$ and  an $\epsilon >0$, 
such that for $0<\gamma <\gamma_0$ and $i=1,2$, if the minimal graph $G(A)$ for the
set $$A=Q\cup\{a_1,e_1(\gamma ),f_1(\gamma ),a_2,e_2(\gamma ),f_2(\gamma )\}$$ is in the
$\epsilon$-neighborhood of the hexagon $H$, then $G(A)$ is symmetric with respect
to the $xz$-plane.
\end{coro}

Finally, in this sequence of lemmas we have:

\begin{lemma}\label{splitfinal}  There is a $\gamma_0 >0$,  such that
for $0<\gamma <\gamma_0$ and $i=1,2$,  the minimal graph $G(A)$ for the set $$A=Q\cup
\{a_1,e_1(\gamma ),f_1(\gamma ),a_2,e_2(\gamma ),f_2(\gamma )\}$$ is unique and
consists of  the two minimal trees $T(\{a_i, b_i, e_i(\gamma),f_i(\gamma )\})$, $i=1,2$.
\end{lemma}

\begin{lemma}\label{circles}  
There is a $\delta >0$ such that if $b'_2\in P$ is a point below the $x$-axis and
${\rm d}(b_2,b'_2)<\delta$, then the length of the minimal tree
$|T(\{a_2,v,c_2(\gamma)\})|$ is  a strictly monotone function of $v\in {\rm
arc}_P(b_2,b'_2)$ and attains its maximum at $b_2$.
\end{lemma}

\begin{figure}[ht]
\centerline{\epsfbox{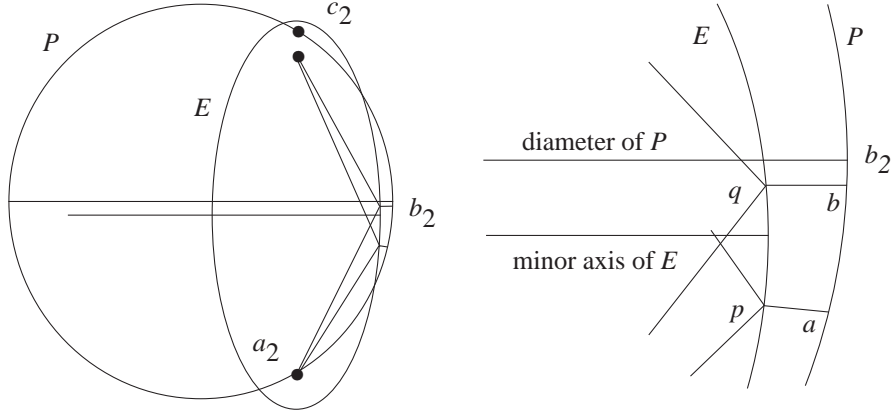}}
\caption{  Perturbation.}
\label{fcircles}
\end{figure}

\begin{proof} $P$ is the circle circumscribed around the hexagon $H$. 
Let $a$ and $b$ be points in   ${\rm arc}_P(b_2,b'_2)$ with $b\in
{\rm arc}_P(b_2,a)$. Note that $T(\{a_2,b,c_2(\gamma)\})$ is a simple triod; denote
the vertex of order 3 by $q$. Let $E$ be the ellipse with foci $a_2$ and
$c_2(\gamma)$, and passing through the point $q$. Let $p$ be the point on $E$ closest
to $a$, see Figure ~\ref{fcircles}. We have

$$
|T(\{a_2,a,c_2(\gamma)\})|<{\rm d}(a,p) + {\rm d}(p,a_2) + {\rm d}(p,c_2(\gamma))=$$
$${\rm d}(a,p) + {\rm d}(q,a_2) + {\rm d}(q,c_2(\gamma))<|T(\{a_2,b,c_2(\gamma)\})|.
$$ 
\end{proof}

\begin{coro}\label{splitcircle}  
There are a $\gamma >0$ and a $\delta >0$ such that if $b'_2\in P$ is a point below the
$x$-axis and ${\rm d}(b_2,b'_2)<\delta$, then the length of the minimal tree
$|T(\{a_2,v,e_2(\gamma),f_2(\gamma)\})|$ is  a strictly monotone function of
$v\in {\rm arc}_P(b_2,b'_2)$ and attains its maximum at $b_2$.
\end{coro}

\section{A knotted minimal tree}

In this section, we use the spherical coordinates, denoted by 
$(r,\theta , \phi)_s$, where  $x=r\sin\phi \cos \theta$, 
$y=r\sin \phi \sin\theta$, and $z=r\cos \phi $. 

Let $d_1(\delta)=(1,\pi,\frac{\pi}{2}-\delta)_s$,
$d_2(\delta)=(1,0,\frac{\pi}{2}+\delta)_s$, and

$$
\begin{array}{rcl}
M_1(\delta)&=&{\rm arc}_{S^2}(d_1,(1,\frac{5\pi}{4},\frac{\pi}{2}-\delta)_s),\\
M_2(\delta)&=&{\rm
arc}_{S^2}((1,\frac{5\pi}{4},\frac{\pi}{2}-\delta)_s,(1,\frac{7\pi}{4},\frac{\pi}{2})_s),\\
M_3&=&Q- {\rm
arc}_{S^2}((1,\frac{5\pi}{4},\frac{\pi}{2})_s,(1,\frac{7\pi}{4},\frac{\pi}{2})_s),\\
M_4(\delta)&=&{\rm
arc}_{S^2}((1,\frac{5\pi}{4},\frac{\pi}{2})_s,(1,\frac{7\pi}{4},\frac{\pi}{2}+\delta)_s),\\
M_5(\delta)&=&{\rm arc}_{S^2}((1,\frac{7\pi}{4},\frac{\pi}{2}+\delta)_s,d_2),\\
M(\delta)&=&\cup_{i=1}^5M_i(\delta).
\end{array}
$$ 

\noindent Thus $M(\delta)$ is a path from $d_1(\delta)$ to $d_2(\delta)$, close
to $Q$, and of length approximately  $\frac{3}{2}$ times the length of $Q$.

\begin{figure}[ht]
\centerline{\epsfbox{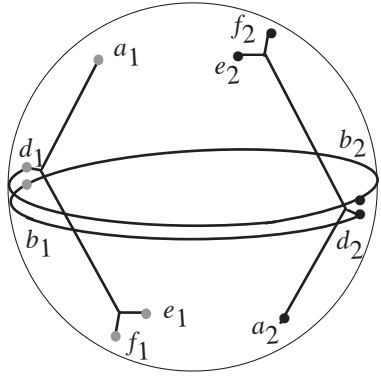}}
\caption{ A design for a knotted minimal tree.}
\label{fknotted}
\end{figure}

\begin{lemma}\label{model} There are a $\gamma $ and a $\delta$  such that the
minimal graph $G(B)$ for the set
$$B=M(\delta)\cup \{a_1,e_1(\gamma),f_1(\gamma),a_2,e_2(\gamma),f_2(\gamma)\}$$ is
unique and consists of the two minimal trees
$T(\{a_i,d_i(\delta),e_i(\gamma),f_i(\gamma)\})$,
$i=1,2$.
\end{lemma}

\begin{proof}  By Lemma ~\ref{splitfinal}, for if  $\gamma $ is sufficiently small,
then the minimal graph $G(A)$ for the set 
$$A=Q\cup\{a_1,e_1(\gamma ),f_1(\gamma ),a_2,e_2(\gamma ),f_2(\gamma )\}$$ 
is unique and 
$$G(A)=T(\{a_1, b_1, e_1(\gamma),f_1(\gamma
)\})\cup T(\{a_2, b_2, e_2(\gamma),f_2(\gamma )\}).$$ Since the sets $A$ and $B$ are
Hausdorff $\delta$-close, there is a $\delta$ such that $$G(B)=T(\{a_1, m_1,
e_1(\gamma),f_1(\gamma )\})\cup T(\{a_2, m_2, e_2(\gamma),f_2(\gamma )\})$$ for some
$m_1$ and $m_2$ on $M$. By Lemma ~\ref{split}, the points $m_i$ are on the $xz$-plane. For
sufficiently small $\delta$, by Corollary ~\ref{splitcircle},  $m_i=d_i(\delta )$.
\end{proof}

For some positive constants $\gamma $ and $\delta$ chosen so that Lemma
~\ref{model} is satisfied,  let $M=M(\delta )$, $d_i=d_i(\gamma )$,  $e_i=e_i(\gamma )$,
and $f_i=f_i(\gamma )$.  For $\epsilon >0$, let $t_2, \dots ,t_{n(\epsilon )-1}$ be points
on $M$ such that
${\rm d}(t_i,t_{i+1}) <\epsilon$ for $i=1,\ldots , n(\epsilon )-2$, 
 $t_2=(1,\pi+\epsilon ,\frac{\pi}{2}-\delta)_s$, and $t_{{n(\epsilon )}-1}=(1,-\epsilon
,\frac{\pi}{2}+\delta)_s$.  In addition let $t_1=(1,\pi-\epsilon ,\frac{\pi}{2}-\delta)_s$
and $t_{n(\epsilon )} =(1,\epsilon ,\frac{\pi}{2}+\delta)_s$.

\begin{thm}\label{example}  
There is an $\epsilon >0$ such that the minimal tree for the set 
$$X=\{a_1,e_1,f_1,a_2,e_2,f_2, t_1,\ldots ,t_{n(\epsilon )}\}$$ 
is unique and knotted.
\end{thm}

\begin{proof} Let $A_1=\{a_1, t_1,t_2, e_1,f_1\}$ and $A_2=\{a_2, t_{n(\epsilon
)-1},t_{n(\epsilon )}, e_2,f_2\}$.  For small $\epsilon$, $T(X)$ contains two subgraphs
close to the minimal  trees $T(\{a_1, d_1, e_1,f_1\})$ and $T(\{a_2, d_2, e_2,f_2\})$.
Since ${\rm d}(t_1,t_2)={\rm d}(t_{n(\epsilon)-1},t_{n(\epsilon)})$ is approximately
$2\epsilon$ and  ${\rm d}(t_i,t_{i+1})<\epsilon$ for the remaining points $t_i$,  
for sufficiently small $\epsilon$, the two subgraphs are $T(A_1)$ and $T(A_2)$. Thus
there is an $\epsilon >0$ such that $T(X)$ is unique and 
$$T(X) =T(A_1)
 \cup [t_2,t_3]\cup \ldots \cup [t_{n(\epsilon )-2},t_{n(\epsilon )-1}]
 \cup T(A_2).
$$ Such $T(X)$ is knotted, see Figure ~\ref{fknotted}.
\end{proof}

\noindent {\bf Remark.} A slight change of the arc $M$ and an appropriate choice of the sequence
of the points $t_i$ can yield a finite set in $S^2$ with two minimal trees, one knotted and
the other unknotted.

\medskip

The example of the knotted minimal tree raises the following questions: 

\begin{enumerate}

\item (M.~Freedman)  What  does the set of $k$-tuples in $S^2$ whose minimal tree is
unknotted look like in the $k$-fold product  $S^2\times\cdots\times S^2$? In particular, what
is the  measure of this set?

\item What is the minimum number $k$ for which there is a set of points in $S^2$ whose 
minimal tree is knotted?

\item (W.~Kuperberg) What is the minimum number of vertices of order 3 in a knotted minimal
tree for a finite subset of $S^2$? The described knotted tree has 6 vertices of order 3.

\item (G.~Kuperberg) There is a finite set whose minimal tree is knotted on the surface of an
ellipsoid with one of the axes much shorter than the other two axes.  What are the
strictly convex closed surfaces in ${\mathbb R}^3$ containing  a finite set whose minimal tree is
knotted? Can any knot be realized in a minimal tree of a finite set on some  convex
surface?

\end{enumerate}

The author would like to thank Greg  Kuperberg for helpful discussions and
William Chen for bringing the problem to her attention.

\end{document}